\documentclass{article}
\usepackage{amssymb,amsmath,multirow}
\usepackage[dvips]{graphicx}
\usepackage{epsfig}
\usepackage{graphics}
\newtheorem{thm}{Theorem}
\newtheorem{lem}{Lemma}
\newtheorem{cor}{Corollary}
\newtheorem{algo}{Algorithm}
\newtheorem{exam}{Example}
\def\pn{\par\smallskip\noindent}
\def\proof{\pn {Proof.} }
\def\endproof{\hfill \quad{$\Box$}\smallskip}

\begin{document}

\title{Double Well Potential Function and Its Optimization in the
 $n$-dimensional Real Space -- Part II}
\author{  {\small Yong Xia }  \\
{\it \small  State Key Laboratory of Software Development
              Environment, }\\
{\it\small School of Mathematics and System Sciences, Beihang University, P. R. China}\\
{ \small Ruey-Lin Sheu}\\
{\it\small Department of Mathematics, National Cheng Kung University, Taiwan}\\
{\small Shu-Cherng Fang}\\
{\it\small  Department of Industrial and Systems Engineering, North Carolina State University, USA}\\
{\small Wenxun Xing}\\
{\it\small  Department of Mathematics Sciences, Tsinghua University, P. R. China}
 }

\date{}

\maketitle
\begin{abstract}
In contrast to taking the dual approach for finding a global minimum solution of
a double well potential function, in Part II of the paper, we
characterize a local minimizer, local maximizer, and
global minimizer directly from the primal side. It is proven
that, for a ``nonsingular" double well function, there exists at most one
local, but non-global, minimizer and at most one local maximizer. Moreover,
when it exists, the local maximizer is ``surrounded" by local
minimizers in the sense that the norm of the local maximizer is
strictly less than that of any local minimizer. We also establish some
necessary and sufficient optimality conditions for the global
minimizer, local non-global minimizer and local
maximizer by studying a convex secular function over specific
intervals. These conditions lead to three algorithms for
identifying different types of critical points of a given
double well function.

%
\end{abstract}

{\bf Keywords}: Double well potential, Local minimizer, Local
maximizer, Global minimum.

\section{Introduction}

In Part I, the double well potential problem (DWP) is defined by
\begin{equation}
\min_{x\in R^n} \left\{
\Pi(x)=\frac{1}{2}\left(\frac{1}{2}\|Bx-c\|^2-d \right)^2 + \frac{1}{2}x^TAx-f^Tx \right\}, \label{d-well}
\end{equation}
where $A$ is an $n\times n$ real symmetric matrix, $B\not= 0$ is an
$m\times n$ real matrix, $c\in R^m$, $d\in R$ and $f\in R^{n}$.
By introducing a continuous variable transformation
$\xi=\frac{1}{2}\|Bx-c\|^2-d$, the double well potential problem (DWP)
can be transformed into the following equivalent
quadratic program over one nonhomogeneous quadratic constraint
(QP1QC):
\begin{eqnarray}
&\min & \Pi(x,\xi)= \frac{1}{2}\xi^2+\frac{1}{2}x^TAx-f^Tx \label{qp1qc:1}\\
&{\rm s.t.}& \xi=\frac{1}{2}\|Bx-c\|^2-d, ~~x \in R^n.\label{qp1qc:2}
\end{eqnarray}
The dual problem of (QP1QC) and the dual of the dual were
studied in Part 1 (Theorem 1) in order to find a global minimum solution to problem (DWP).

For practical applications, knowing only the global minimum of a double well
potential function may not be sufficient. For example, the double
well potential model can be used to describe the ion-molecule
reactions, where the intermediate molecule complexes must go across the
energy barrier to cause reactions \cite{Brauman}. Researchers have to know the potential
difference between the energy wells (caused by local minima) and energy
barrier (caused by local maximum). The understanding of all types of critical
points of a double well function is
thus necessary.

Mathematically, we are motivated by the pioneering work of Mart\'{\i}nez \cite{Mar}
which showed that a trust-region subproblem (TRS) \cite{Conn} of the following form
\begin{eqnarray}
&\min & \frac{1}{2}x^TAx-f^Tx \label{tr:1}\\
&{\rm s.t.}& \|x\|^2 =\Delta, ~~x \in R^n \label{tr:2}
\end{eqnarray}
(with $\Delta$ being a positive scalar)
has at most one local, but non-global, minimizer. Please notice that, on one hand,
problem (QP1QC) can be regarded as an extension
of problem (TRS) towards the nonhomogeneous and possibly
singular case. On the other hand, the penalty version of the trust-region
subproblem, namely,
$$
\min_{x \in R^n}  \frac{1}{2}x^TAx-f^Tx+\theta(\|x\|^2-\Delta)^2
$$
(with the penalty parameter $\theta$ being sufficiently large) is clearly a
special case of the double well potential problem (DWP). Therefore, our
approach to analyzing the local non-global minimizer of a double well
potential problem extends the results of \cite{Mar}. Moreover, when restricted to problem (TRS), our approach simplifies the proof provided in \cite{Mar}. Although, in general, a double well potential
problem may have infinitely many local, but non-global, minimizers (see
Figure 1), we'll show that, after taking the space reduction technique developed in Section 2 of Part I,
the reduced nonsingular problem has at most one local non-global
minimizer and at most one local maximizer.

\begin{figure}[!hbpt]\label{fig:0}
\begin{center}
\includegraphics[bb=0 0 544 496,width=0.45\textwidth]{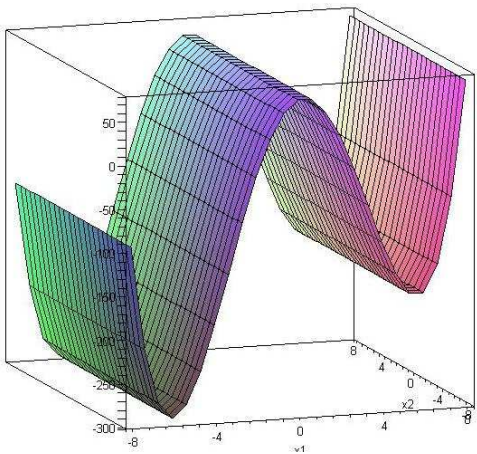}
\end{center}
\caption{A double well potential problem having infinitely many local non-global minima.}
\end{figure}



We remark that characterizing the local maximizer of the
trust-region subproblem (\ref{tr:1})-(\ref{tr:2}) can be reduced to the problem of
finding a local minimizer of
(\ref{tr:1}) with $A$ being replaced by $-A$. However, due to the non-symmetric nature, it is no
longer the case for the double well potential problem
(\ref{d-well}). Hence Mart\'{\i}nez's approach may not be able to characterize the
local maximizer for a general (DWP) problem.

%
%
%


In the rest of the paper, a characterization of the local, but
non-global, minimizer of a double well function is provided in
Section 2. Then, a characterization of the global minimizer of a
double well function is given in Section 3, while the local
maximizer is characterized in Section 4. Computational algorithms
for each type of the optimizers of a double well potential function
are proposed in Section 5 with some illustrative examples. Some
concluding remarks are given in Section 6.

Here we define some notations to be used throughout the paper.
Let $S^n$ be the set of all $n$-dimensional symmetric real matrices,
$S_+^n$ be the set of all $n$-dimensional positive semi-definite matrices,
and $S_{++}^n$ be the set of all n-dimensional positive definite matrices.
For any $P, Q \in S^n$, $P\succeq Q$ means that matrix $P-Q \in S_+^n$
and  $P\succ Q$ means that matrix $P-Q \in S_{++}^n$.
We sometimes write $Q\preceq P$ for $P\succeq Q$
and $Q\prec P$ for $P\succ Q$. The $i$th smallest eigenvalue of $P \in S^n$ is
denoted by $\sigma_i(P)$ and the determinant of $P$ by $\det(P)$.
The $n$-dimensional identity matrix is denoted by $I$. For a vector $x\in R^n$,
$x\geq 0$ ($x > 0$) means that each component of $x$ is nonnegative (positive) and Diag$(x)$
is an $n$-dimensional diagonal matrix with diagonal components being
$x_1,\ldots,x_n$. Moreover, for a number $\beta \in R^n$, sign$(\beta)=\frac{\beta}{|\beta|}$ if $\beta \neq 0$, otherwise
sign$(\beta)=0$.

\section{Characterization of local non-global minimizer}

Following the space reduction technique developed in Part I, without
loss of generality, we may assume that, in problem (BWP), $B^TB$ is positive definite such that,
matrices $A$ and $B^TB$ are simultaneously diagonalizable
via congruence, i.e., there is a nonsingular matrix $P$ such that
$D\triangleq P^TAP={\rm Diag} (\alpha_1,\ldots,\alpha_n)$ with
$\alpha_1\leq\ldots\leq\alpha_n$ and $P^TB^TBP= I$. It follows
immediately that $(B^TB)^{-1}=PP^T$. Let
\[
w=P^{-1}x-P^TB^Tc, \]
then we have
\begin{eqnarray*}
\frac{1}{2}\|Bx-c\|^2-d &=&\frac{1}{2}\|B(Pw+PP^TB^Tc)-c\|^2-d\\&=&
\frac{1}{2}w^Tw+\frac{1}{2}c^T(I-B(B^TB)^{-1}B^T)c-d
\end{eqnarray*}
and
\[
\frac{1}{2}x^T Ax-f^Tx =\frac{1}{2}w^TDw+c^TBPDw+\frac{1}{2}c^TBPDP^TB^Tc-f^TPw-f^TPP^TB^Tc.
\]
For simplicity, we define \
$\nu=-\frac{1}{2}c^T(I-B(B^TB)^{-1}B^T)c+d$ \ and \
$\psi=P^Tf-DP^TB^Tc$. By dropping the constant terms, we can rewrite problem
(DWP) defined in (\ref{d-well}) as

\begin{equation}
\min \left\{g(w)=
\frac{1}{2}\left(\frac{1}{2}\|w\|^2-\nu\right)^2+\frac{1}{2}w^TDw-\psi^Tw\right\}.
\label{d-well:4}
\end{equation}

Recall that the canonical primal problem defined in (19) of Part I is to
minimize
\begin{equation}
\frac{1}{2}\left(\frac{1}{2}\|w\|^2-\varphi^Tw-\nu\right)^2+\frac{1}{2}w^TDw-\psi^Tw.
\label{d-well:19}
\end{equation}
The form in (\ref{d-well:4}) is a further simplified version of form
(\ref{d-well:19}) by setting $\varphi=0$. In this way, the third
order term in problem (DWP) is eliminated and the complexity is decreased for analysis.
It's interesting to note that, in the
finite deformation theory, the diagonal matrix $D$ represents the
material constants, the first order coefficient vector $\psi$ stands
for the external forces, and the Cauchy-Green strain
$\frac{1}{2}\|w\|^2-\nu$ measures the square of local changes in
distance due to deformation. As we shall observe below, the first
order and the second order necessary conditions of (\ref{d-well:4})
(see \cite{NW}) are highly related to the term of
$(\frac{1}{2}\Vert {w}\Vert^2-\nu)I+D$, which is the sum of the
Cauchy-Green strain and the material constants. Our first result of
Lemma \ref{lem:2} will show that, at a local minimum of the double well
potential function, the Cauchy-Green strain can not be too small, at
least no smaller than the negative of the second smallest material
constant.

\begin{lem} \label{lem:1}
Assume that $\underline{w}$ is a local minimizer of
(\ref{d-well:4}). It holds that
\begin{eqnarray}
&&\nabla g(\underline{w})=\left((\frac{1}{2}\Vert
\underline{w}\Vert^2-\nu)I+
D\right)\underline{w}-\psi=0\label{con:1},\\
&&\nabla^2 g(\underline{w})=\underline{w}\underline{w}^{T}
+(\frac{1}{2}\Vert \underline{w}\Vert^2-\nu)I+D\succeq
0.\label{con:2}
\end{eqnarray}
\end{lem}

\begin{lem} \label{lem:2}
Assume that $n\geq 2$ and $\underline{w}$ is a local minimizer of
(\ref{d-well:4}). It holds that
\begin{eqnarray}
\frac{1}{2}\Vert \underline{w}\Vert^2-\nu+\alpha_2\geq 0.
\label{lem:eq}
\end{eqnarray}
Furthermore, if $\alpha_1<\alpha_2$, then
\begin{eqnarray} \label{lem:ineq}
\frac{1}{2}\Vert \underline{w}\Vert^2-\nu+\alpha_2> 0.
\end{eqnarray}
\end{lem}
\proof Suppose that the statement (\ref{lem:eq}) is false, then $\frac{1}{2}\Vert
\underline{w}\Vert^2-\nu+\alpha_2<0$. Hence $\frac{1}{2}\Vert
\underline{w}\Vert^2-\nu+\alpha_1<0$. Let $e_1^T=(1,0,0,\ldots,0)$
and $e_2^T=(0,1,0,\ldots,0)$. If
$e_1^T\underline{w}=\underline{w}_1=0$, by the necessary condition
(\ref{con:2}), we have
\begin{eqnarray}
0\leq e_1^T(\underline{w}\underline{w}^{T}+(\frac{1}{2}\Vert
\underline{w}\Vert^2-\nu)I+D)e_1
=\frac{1}{2}\Vert\underline{w}\Vert^2-\nu+\alpha_1<0,
\end{eqnarray}
which causes a contradiction. On the other hand, if $e_1^T\underline{w}\neq 0$, then, by
(\ref{con:2}) again, we have
\begin{eqnarray*}
0&\leq
&((-\underline{w}_2)e_1+(\underline{w}_1)e_2)^T(\underline{w}\underline{w}^{T}
+(\frac{1}{2}\Vert \underline{w}\Vert^2-\nu)I+D)((-\underline{w}_2)e_1+(\underline{w}_1)e_2)\\
&=&(\frac{1}{2}\Vert\underline{w}\Vert^2-\nu+\alpha_1)(\underline{w}_2)^2+(\frac{1}{2}\Vert
\underline{w}\Vert^2-\nu+\alpha_2)(\underline{w}_1)^2<0.
\end{eqnarray*}
It again causes a contradiction. Therefore, the statement (\ref{lem:eq}) must be true.

When $\alpha_1<\alpha_2$, suppose that the statement (\ref{lem:ineq}) is false, then we have
\begin{eqnarray}\label{lem:eq_1}
\frac{1}{2}\Vert \underline{w}\Vert^2-\nu+\alpha_2=0.
\end{eqnarray}
By (\ref{lem:eq_1}), we know that the second order necessary condition
(\ref{con:2}) becomes
\begin{eqnarray}\label{x_2}
0\preceq \underline{w}\underline{w}^{T}+(\frac{1}{2}\Vert
\underline{w}\Vert^2-\nu)I+D = \underline{w}\underline{w}^{T}+D-
\alpha_2I.
\end{eqnarray}
Since the first two leading principal minors of the matrix in (\ref{x_2}) are
nonnegative, we have
\begin{eqnarray}\label{x_0}
 \underline{w}_1^{2}+\alpha_1-\alpha_2\ge0
\end{eqnarray}
and
\begin{equation}\label{x1}
\det\left\{\left[\begin{array}{cc} \underline{w}_1^{2}&
\underline{w}_1\underline{w}_2
\\ \underline{w}_1\underline{w}_2& \underline{w}_2^{2}
\end{array}\right]+\left[\begin{array}{cc}\alpha_1-\alpha_2&0\\0&0 \end{array}\right]
\right\}= (\alpha_1-\alpha_2)\underline{w}_2^{2}\ge 0.
\end{equation}
Remember that $\alpha_1-\alpha_2<0$, inequality (\ref{x_0}) implies that
$\underline{w}_1\not=0$. Moreover, inequality (\ref{x1}) implies that
$\underline{w}_2=0$. Together with (\ref{con:1}), we obtain that
$\psi_2=0$ and
$$\underline{w}_1=\frac{2\psi_1}{\Vert \underline{w}\Vert^2-2\nu+2\alpha_1}
=\frac{-\psi_1}{\alpha_2-\alpha_1}.
$$
Without loss of generality, we assume that $\psi_1<0$, and hence
$\underline{w}_1>0$. This implies, from (\ref{lem:eq_1}) and the fact that
$\underline{w}_2=0$, we have
$$\underline{w}_1=\sqrt{2\nu-2\alpha_2-\sum_{i=3}^n\underline{w}_i^{2}}.$$
Consider the following parametric curve in $\mathbb{R}^n$:
\begin{equation}\label{curve}
\gamma(t)=\{(k(t),t,\underline{w}_3,\ldots,\underline{w}_n)\vert
k(t)=\sqrt{
2\nu-2\alpha_2-t^2-\sum_{i=3}^n\underline{w}_i^{2}}=\sqrt{\underline{w}_1^{2}-t^2},
~t\in R\}
\end{equation}
where $\gamma(0)=\gamma(\underline{w}_2)=\underline{w}$, i.e.,
$\gamma(t)$ passes through $\underline{w}$ at $t=0$. Evaluating $g(w)$ on
$\gamma(t)$, we have
\begin{eqnarray*}
g(\gamma(t))&=&\frac{1}{2}\left(\frac{k(t)^2}{2}+\frac{t^2}{2}
+\frac{1}{2}\sum_{i=3}^n\underline{w}_i^{2}-
\nu+\alpha_1
\right)^2-\alpha_1\left(\frac{k(t)^2}{2}+\frac{t^2}{2}+
\frac{1}{2}\sum_{i=3}^n\underline{w}_i^{2}-\nu\right)-\frac{\alpha_1^2}{2}\\
& &+\frac{1}{2}\left(\alpha_1k(t)^2+\alpha_2t^2+\sum\limits_{i=3}^n\alpha_i
\underline{w}_i^{2}\right)
-\psi_1k(t)-\sum_{i=3}^n\psi_i\underline{w}_i\\
&=&\frac{(\alpha_2-\alpha_1)^2}{2}+\frac{\alpha_2-\alpha_1}{2}t^2+
\sum_{i=3}^n\frac{\alpha_i-\alpha_1}{2}\underline{w}_i^{2}-\psi_1\sqrt{\underline{w}_1^{2}-
t^2}
-\sum_{i=3}^n\psi_i\underline{w}_i+\alpha_1\alpha_2-\frac{\alpha_1^2}{2}.
\end{eqnarray*}
It is not difficult to see that $t=0$ is a local minimum point of $g(\gamma(w))$
since $\underline{w}$ is a local minimizer of $g(w)$. However, this conclusion
contradicts to the fact that
$$\frac{d}{dt}g(\gamma(0))=\frac{d^2}{dt^2}g(\gamma(0))=\frac{d^3}{dt^3}g(\gamma(0))=0;
~\frac{d^4}{dt^4}g(\gamma(0))=-\frac{3(\alpha_2-\alpha_1)}{\underline{w}_1^{2}}<0.$$
Therefore, the statement (\ref{lem:ineq}) must be true, if
$\alpha_1<\alpha_2$.
%
\endproof

The next result Lemma shows that any critical point of
the double well potential function having a sufficiently large Cauchy-Green
strain (larger than the negative of all the material constants) must
be a global minimum point.

\begin{lem} \label{lem:3}
Let $w^*$ be a critical point of the function $g(w)$ in problem (\ref{d-well:4}) with $\nabla
g(w^*)=0$. If
\begin{eqnarray}
\frac{1}{2}\Vert w^{*}\Vert^2-\nu+\alpha_1 \geq 0,\label{con:3}
\end{eqnarray}
then $w^*$ is a global minimizer of problem (\ref{d-well:4}). In particular,
a local minimizer $\underline{w}$ of problem (\ref{d-well:4}) satisfying condition
(\ref{con:3}) must be a global minimizer.
\end{lem}
\proof Define $Q=(\frac{1}{2}\Vert w^{*}\Vert^2-\nu)I+D$. By the
assumption that $\frac{1}{2}\Vert w^{*}\Vert^2-\nu+\alpha_1 \geq 0$,
it follows that $\frac{1}{2}\Vert w^{*}\Vert^2-\nu+\alpha_i \geq 0,\
\forall i\in[1:n]$, and $Q$ is positive semidefinite. Then,
\begin{eqnarray}
 g(w)&=&  \frac{1}{2}\left(\frac{1}{2}\|w\|^2-\nu\right)^2+\frac{1}{2}w^TDw-\psi^Tw \nonumber \\
 &=& \frac{1}{2}\left(\frac{1}{2}\|w\|^2-\nu\right)^2-
 \frac{1}{2}w^T(\frac{1}{2}\Vert w^{*}\Vert^2-\nu)Iw+\frac{1}{2}w^TQw-\psi^Tw\nonumber\\
 &=& \frac{1}{8}\Vert w\Vert^4-\frac{1}{4}\|w\|^2\Vert w^{*}\Vert^2
 +\frac{\nu^2}{2} +\frac{1}{2}w^TQw-\psi^Tw\nonumber\\
 &=&\frac{1}{8}\left(\|w\|^2-\Vert w^{*} \Vert^2\right)^2+ \frac{1}{2}w^TQw-\psi^Tw
 -\frac{1}{8}\Vert w^{*}\Vert^4 +\frac{\nu^2}{2} \label{fun:0}\\
 &\geq& \frac{1}{2}w^TQw-\psi^Tw-\frac{1}{8}\Vert w^{*}\Vert^4 +\frac{\nu^2}{2}\label{fun:1}\\
&\geq& \frac{1}{2}w^{*T}Qw^*-\psi^Tw^*-\frac{1}{8}\Vert w^{*}\Vert^4 +\frac{\nu^2}{2}\label{fun:2}\\
&=& g(w^*).\label{fun:3}
\end{eqnarray}
Since $Q\succeq0$, the lower bound function expressed in (\ref{fun:1}) is a convex
quadratic function. Its global minimum is attained at any $\hat w$ satisfying
$Q{\hat w}-\psi=(\frac{1}{2}\Vert w^{*}\Vert^2-\nu)I{\hat w}+D{\hat
w}-\psi=0$. Since $w^*$ is a critical point of (\ref{d-well:4}), by equation (\ref{con:1}) in
Lemma \ref{lem:1}, it is a global minimizer
of the lower bound function in (\ref{fun:1}) and thus inequality (\ref{fun:2}) holds.
Finally, (\ref{fun:3}) becomes true by substituting $w^*$ into
(\ref{fun:0}).
%
%
\endproof

\begin{thm}\label{theorem1}
The double well potential problem (\ref{d-well:4}) has at most one
local, but non-global, minimizer.
\end{thm}
\proof Let us assume that $n\geq 2$ first. Lemmas \ref{lem:2} and
\ref{lem:3} imply that any local, but non-global, minimizer $\underline{w}$
of problem (\ref{d-well:4}) exits only if
$\alpha_1<\alpha_2$ and $-\alpha_2< \frac{1}{2}\Vert
\underline{w}\Vert^2-\nu<-\alpha_1$. Consequently, we know the matrix
$(\frac{1}{2}\Vert \underline{w}\Vert^2-\nu)I+D$ is nonsingular with
its first diagonal element being negative and others positive.
Therefore, $\underline{w}$ can be uniquely determined by equation
(\ref{con:1}) with
\begin{equation}
\underline{w}_i=\frac{2\psi_i}{\Vert
\underline{w}\Vert^2-2\nu+2\alpha_i}, ~~i\in[1:n].\label{sol}
\end{equation}
From (\ref{con:2}), we have $2\underline{w}^{2}_1+ \Vert
\underline{w}\Vert^2-2\nu+2\alpha_1\geq 0.$ Since $\Vert
\underline{w}\Vert^2-2\nu+2\alpha_1<0$, we know that
\begin{equation}\label{a1}
\underline{w}_1\not=0,~ \psi_1\not=0~ \hbox{  and  }
2\nu-2\alpha_1>\Vert \underline{w}\Vert^2>0.
\end{equation}
Putting all $\underline{w}_i$ together, we have
 $$
\sum_{i=1}^n\frac{4\psi_i^2}{(\Vert
\underline{w}\Vert^2-2\nu+2\alpha_i)^2}=\Vert \underline{w}\Vert^2.
 $$
In other words, the norm square of the local minimizer, i.e., $\Vert
\underline{w}\Vert^2$, must be the root of the following secular
function on a specific open interval:
\begin{equation}
 h(t)=\sum_{i=1}^n\frac{4\psi_i^2}{(t-2\nu+2\alpha_i)^2}-t,
 ~~t\in(\max\{2\nu-2\alpha_2,0\},2\nu-2\alpha_1).\label{varphi}
\end{equation}
Notice that each root of $h(t)=0$ can only correspond to one local non-global
minimizer of problem (\ref{d-well:4}) using (\ref{sol}). Taking a simple
calculation of (\ref{varphi}), we have
\begin{eqnarray}
h'(t)&=&-\sum_{i=1}^n\frac{8\psi_i^2}{(t-2\nu+2\alpha_i)^3}-1,\label{a4}\\
h''(t)&=&\sum_{i=1}^n\frac{24\psi_i^2}{(t-2\nu+2\alpha_i)^4}>0.\label{a3}
\end{eqnarray}
Therefore, the secular function $h(t)$ is a strictly convex
function on $(\max\{2\nu-2\alpha_2,0\},2\nu-2\alpha_1)$ with at most
two roots. Furthermore, since $(\Vert
\underline{w}\Vert^2-2\nu)I+2D$ is nonsingular, the second order
necessary condition (\ref{con:2}) implies that
$$
2(\Gamma\underline{w})(\Gamma\underline{w})^T+ {\rm Diag}(-1,1,\ldots,1)  \succeq 0,
$$
where
\begin{equation}
\Gamma={\rm Diag}\left(\frac{1}{\sqrt{-\Vert
\underline{w}\Vert^2+2\nu-2\alpha_1}}, \frac{1}{\sqrt{\Vert
\underline{w}\Vert^2-2\nu+2\alpha_2}},\ldots, \frac{1}{\sqrt{\Vert
\underline{w}\Vert^2-2\nu+2\alpha_n}} \right).\label{gamma:1}
\end{equation}
Since its determinant is nonnegative, we have
\begin{eqnarray}
0&\leq& {\rm det}\left(2(\Gamma \underline{w}) (\Gamma \underline{w})^T+
{\rm Diag}(-1,1,\ldots,1)\right)\label{determinant1} \\
&=&{\rm det}({\rm Diag}(-1,1,\ldots,1))\cdot{\rm det}\left( 2{\rm Diag}(-1,1,\ldots,1)
(\Gamma \underline{w}) (\Gamma \underline{w})^T+I\right)\nonumber \\
&=&-1\cdot(2(\Gamma\underline{w})^T{\rm Diag}(-1,1,\ldots,1)(\Gamma\underline{w})+1)\nonumber \\
&=&-\sum_{i=1}^n\frac{8\psi_i^2}{(\Vert
\underline{w}\Vert^2-2\nu+2\alpha_i)^3}-1\nonumber\\
&=&h'(t)\mid_{t=\Vert \underline{w}\Vert^2}.\label{determinant2}
\end{eqnarray}
In other words, if $\underline{w}$ is a local minimizer of problem (\ref{d-well:4}), then it
must satisfy the second order necessary condition (in matrix form)
whose determinant is the first derivative of the secular function at
$\Vert\underline{w}\Vert^2$. However, a strictly convex
function has at most one root with a nonnegative first
derivative. Thus we have shown the theorem for $n\ge2$. When
$n=1$, it amounts to setting $\alpha_2=\infty$ in the above
analysis, and the proof follows.
\endproof

The next corollary provides some simple sufficient conditions for having no
local non-global minimizer.
\begin{cor}\label{cor:1}
When one of the following conditions is met:
\begin{itemize}
\item[(i)]
$2\nu-2\alpha_1\leq 0$ (in this case $g(w)$ is convex);
\item[(ii)] $\alpha_1=\alpha_2$;
\item[(iii)] $\psi_1=0$;
\item[(iv)] $\max\{2\nu-2\alpha_2,0\}< 2\nu-2\alpha_1$,
$\psi_1\neq 0$ and
$\min_{t\in(\max\{2\nu-2\alpha_2,0\},2\nu-2\alpha_1)}h(t)>0$;
\end{itemize}
any local minimizer of the double well potential problem
(\ref{d-well:4}) is globally optimal.
\end{cor}
\proof (i) If $2\nu-2\alpha_1\leq 0$, then $\Vert
w\Vert^2-2\nu+2\alpha_i \geq 0$ for any $w\in R^n$ and
$i\in[1:n]$. Using the second derivative of $g(w)$ in
(\ref{con:2}), we have
$$
\nabla^2 g(w)=2ww^{T}+(\Vert w\Vert^2-2\nu)I+2D\succeq 0,
$$
which shows that $g(w)$ is indeed convex and any local optimum becomes globally optimal.

(ii) If
$\underline{w}$ is a local minimizer and
$\alpha_1=\alpha_2$, then inequality (\ref{lem:eq}) in Lemma \ref{lem:2} and Lemma \ref{lem:3} imply that
$\underline{w}$ is indeed the global minimizer of problem (\ref{d-well:4}).

(iii) If $\psi_1=0$, then equation (\ref{con:1}) leads to either $\underline{w}_1=0$
or $\Vert \underline{w}\Vert^2-2\nu+2\alpha_1=0$. By property
(\ref{con:2}), $\underline{w}_1=0$ further implies that $\Vert
\underline{w}\Vert^2-2\nu+2\alpha_1\ge0$. Using Lemma \ref{lem:3}, both cases lead
$\underline{w}$ to be a global minimizer.

(iv) In this case, the secular function $h(t)$ actually does not
have any solution in its domain.
\endproof

The key result of establishing a necessary and sufficient condition for
local, non-global minimizer is provided below.
\begin{thm}\label{cor:0}
The double well potential problem (\ref{d-well:4}) has a
local-nonglobal minimizer if and only if there is a
$\underline{t}^*\in(\max\{2\nu-2\alpha_2,0\},2\nu-2\alpha_1)$ such
that the secular function  $h(\underline{t}^*)=0$ and $
h'(\underline{t}^*)>0$. Moreover, when it exists, the local non-global minimizer is given by
\begin{equation}
\underline{w}
=\left(\frac{2\psi_1}{\underline{t}^*-2\nu+2\alpha_1},\ldots,\frac{2\psi_n}
{\underline{t}^*-2\nu+2\alpha_n} \right). \label{tw}
\end{equation}
\end{thm}
\proof Suppose that
$h(\underline{t}^*)=\sum_{i=1}^n\frac{4\psi_i^2}{(\underline{t}^*-2\nu+2\alpha_i)^2}
-\underline{t}^*=0$ with
$\underline{t}^*\in(\max\{2\nu-2\alpha_2,0\},2\nu-2\alpha_1)$.
For the $\underline{w}$ defined by (\ref{tw}), we have
$\underline{t}^*=\Vert\underline{w}\Vert^2$ and $\underline{w}$
satisfies the first order necessary condition (\ref{con:1}).
Moreover, the diagonal matrix $(\frac{1}{2}\Vert
\underline{w}\Vert^2-\nu)I+D$ is nonsingular with positive diagonal
elements  except for the first one. By Weyl's inequality (see
\cite{Horn}, Theorem 4.3.1), we can estimate the largest $n-1$
eigenvalues of the second order matrix
$\underline{w}\underline{w}^{T}+(\frac{1}{2}\Vert
\underline{w}\Vert^2-\nu)I+D$ by
\begin{eqnarray}
\sigma_i\left(\underline{w}\underline{w}^{T}+(\frac{1}{2}\Vert
\underline{w}\Vert^2-\nu)I+D\right)&\ge&
\sigma_1\left(\underline{w}\underline{w}^{T}\right)+\sigma_i\left(\frac{1}{2}\Vert
\underline{w}\Vert^2-\nu)I+D\right)\nonumber\\
&\ge&\sigma_i\left(\frac{1}{2}\Vert
\underline{w}\Vert^2-\nu)I+D\right)\nonumber\\
&>&0,~~~~~\hbox{for}~i=2,3,\ldots,n.\label{new:1}
\end{eqnarray}
Since $h'(\underline{t}^*)>0$, by (\ref{determinant1}) and
(\ref{determinant2}), we have
\begin{equation}
\det\left(\underline{w}\underline{w}^{T}+(\frac{1}{2}\Vert
\underline{w}\Vert^2-\nu)I+D\right)=\frac{h'(\underline{t}^*)}{2{\det}^2(\Gamma)}>0,
\label{new:2}
\end{equation}
where $\Gamma$ is defined in (\ref{gamma:1}).
Combining (\ref{new:1}) with (\ref{new:2}), we know that the
smallest eigenvalue of the second order matrix must be positive, i.e.,
$\underline{w}\underline{w}^{T}+(\frac{1}{2}\Vert
\underline{w}\Vert^2-\nu)I+D\succ0.$ This is a sufficient
condition to guarantee that $\underline{w}$ is a local minimizer of problem (\ref{d-well:4}).

On the hand side, let $\underline{w}$ be  a local, non-global minimizer of
problem (\ref{d-well:4}), which is unique quaranteed
by Theorem \ref{theorem1}.
Let $\underline{t}^*=\Vert\underline{w}\Vert^2$. By the
proof of Theorem \ref{theorem1},
we know
$\underline{t}^*\in(\max\{2\nu-2\alpha_2,0\},2\nu-2\alpha_1)$. Moreover,
$\underline{w}$ can be expressed by $\underline{t}^*$ as in
(\ref{tw}) because $\underline{w}$ satisfies the first order
necessary condition (\ref{con:1}). Also we have $h(\underline{t}^*)=0$
and $ h'(\underline{t}^*)\ge0$.
It remains for us to show that $ h'(\underline{t}^*)>0$. Suppose that, by
contradiction, $h'(\underline{t}^*)=0$. From (\ref{new:2}), we have
\begin{equation*}
\det\left(\underline{w}\underline{w}^{T}+(\frac{1}{2}\Vert
\underline{w}\Vert^2-\nu)I+D\right)=\frac{h'(\underline{t}^*)}{2{\det}^2(\Gamma)}=0
\end{equation*}
and thus there is a $u=(u_1,\ldots,u_n)^T\neq 0$ such that
\begin{equation}
\underline{w}\underline{w}^{T}u+\left((\frac{1}{2}\Vert
\underline{w}\Vert^2-\nu)I+D \right)u=0.\label{new:000}
\end{equation}
From (\ref{new:000}), we can write
\begin{equation*}
u_i=\frac{-\underline{w}_i(u^T\underline{w})}{\frac{1}{2}\Vert
\underline{w}\Vert^2-\nu+\alpha_i},~i=1,2,\ldots,n.
\end{equation*}
Then, $u\neq 0$ implies that
\begin{equation}
u^T\underline{w}\neq 0.\label{nzero}
\end{equation}
Consider the double well potential function along the direction $u$
defined by $q(\beta):=g(\underline{w}+\beta u)$.
It is routine to verify that
\begin{eqnarray*}
q'(\beta)&=& \nabla g(\underline{w}+\beta u)u, \\
q''(\beta)&=& u^T \nabla^2g(\underline{w}+\beta u)u,\\
q'''(\beta)&=&3u^T (\underline{w}+\beta u )u^Tu  .
\end{eqnarray*}
By the first order necessary condition (\ref{con:1}), we have
$q'(0)=0$. By (\ref{con:2}) and (\ref{new:000}), we further have $q''(0)=0$.
However, (\ref{nzero}) implies that
\[
q'''(0)=3(u^T \underline{w})(u^Tu) \neq 0.
\]
This result contradicts to the fact that $\underline{w}$ is a local minimizer
of problem (\ref{d-well:4}). Therefore, $ h'(\underline{t}^*)>0$ and the proof is complete.
\endproof

\section{Characterization of global minimizer}

In this section, we try to characterize different aspects of the global
minimizer of the double well potential problem.
We first observe that the double well potential function
tends to $+\infty$ as $\Vert w\Vert^2\rightarrow \infty$. Therefore,
the global minimizer of problem (\ref{d-well:4}) always exists. Our first
result is that each component of the global minimizer must be
of the same sign as the corresponding component of the external
force (i.e., the first-order term vector).

\begin{lem}\label{lem:4}
If $w^{*}$ is the global minimizer of (\ref{d-well:4}), then
\begin{equation}
\psi_iw_i^{*}\geq 0,~i\in[1:n].
\end{equation}
\end{lem}
\proof Let $\widetilde{w}=(-w_1^{*},
w_2^{*},w_3^{*},\ldots,w^{*}_n)$. Since the only odd-order term in
$g(w)$ is the linear term, we have
$$
g(w^{*})-g(\widetilde{w})=-\psi_1(w_1^{*}-\widetilde{w}_1)=-2\psi_1w_1^{*}\le0.
$$
Hence we know $\psi_1w^{*}_1\ge0$. A similar argument applies for any other components.
\endproof

The next result shows that the sufficient condition
$\frac{1}{2}\Vert w^{*}\Vert^2-\nu+\alpha_1 \geq 0$ in Lemma
\ref{lem:3} is indeed necessary for a critical point to become the global minimizer.

\begin{thm}\label{thm:suf-nec}
$w^{*}$ is a global minimizer of (\ref{d-well:4}) if and only if
\begin{eqnarray}
\nabla g(w^{*})=\left((\frac{1}{2}\Vert
w^{*}\Vert^2-\nu)+D\right)w^{*}-\psi&=&0\label{glob:s2}
\end{eqnarray}
and
\begin{eqnarray}
\Vert w^{*}\Vert^2-2\nu+2\alpha_1&\geq &0.\label{glob:s1}
\end{eqnarray}
\end{thm}
\proof The sufficiency is clear from Lemma \ref{lem:3}. In addition,
we can observe that the necessity of (\ref{glob:s2}) follows
immediately from equation (\ref{con:1}). It remains to show that
(\ref{glob:s1}) is also a necessary condition.

To avoid triviality, we may assume that $\alpha_1<\alpha_2$.
Otherwise, by substituting $\alpha_1=\alpha_2$ into (\ref{lem:eq}),
we can obtain the result at once. Suppose that $\Vert
w^{*}\Vert^2-2\nu+2\alpha_1 < 0,$ then (\ref{con:2}) implies that
$2{w_1^*}^{2}+ \Vert {w^*}\Vert^2-2\nu+2\alpha_1\geq 0.$ Hence we have
$w_1^{*2}\not=0$. Using (\ref{glob:s2}), we have
$$
2\psi_1w_1^{*} = (\Vert w^{*}\Vert^2-2\nu+2\alpha_1 )w_1^{*2}<0.
$$
This causes a contradiction to Lemma \ref{lem:4} and the proof follows.
\endproof
%

An immediate consequence of Theorem \ref{thm:suf-nec} is that the
sign of the first component of the local non-global minimizer, if it exists,
must be opposite to that of the first component of a global
minimum solution.

\begin{cor}\label{cor:3}
If $\underline{w}$ be the local non-global minimizer and $w^{*}$ is
a global minimizer of $g(w)$ of problem (\ref{d-well:4}), then
\begin{equation}
{\rm sign}(\psi_1)={\rm sign}(w_1^{*})=-{\rm sign}(\underline{w}_1)\in \{-1,~1\}.
\end{equation}
\end{cor}
\proof Since both $\underline{w}$ and $w^{*}$ are critical points,
Theorem \ref{thm:suf-nec} implies that $\Vert
\underline{w}\Vert^2-2\nu+2\alpha_1< 0$ and $\Vert
w^{*}\Vert^2-2\nu+2\alpha_1\geq0$. It follows from condition (iii) of
Corollary \ref{cor:1} that $\psi_1\not=0$ and, from (\ref{con:1}),
\begin{eqnarray}
(\Vert \underline{w}\Vert^2-2\nu+2\alpha_1)\underline{w}_1&=&2\psi_1,\label{sign-3}\\
(\Vert w^{*}\Vert^2-2\nu+2\alpha_1)w_1^{*}&=&2\psi_1.\label{sign-4}
\end{eqnarray}
Consequently, $ {\rm sign}(\underline{w}_1)=-{\rm sign}(\psi_1)=-{\rm
sign}(w_1^{*})\in\{-1,~1\}. $
\endproof
%
%
%

In Section 4 of Part I, we have shown that the dual of the dual of
the canonical primal problem (P) (see equation (19) of Part I) is
equivalent to only a portion of (P) subject to $n$ linear
constraints (see equation (35) of Part I). Moreover, that portion
contains the global minimizer. In the simplified version here, we
have the third order term coefficient
$\varphi=(\varphi_1,\varphi_2,\ldots,\varphi_n)^T=0$, which reduces
the dual of the dual problem in Part I to the following problem:
\begin{eqnarray}\label{dualofdual}
\begin{array}{rll}
P_0^{dd}=&
\inf\limits_{\lambda \in R^n}&P^{dd}(\lambda)=\sum\limits_{i=1}^n\alpha_i\lambda_i
-\sum\limits_{i=1}^n|\psi_i|\sqrt{2\lambda_i}
+\frac{1}{2}(\sum\limits_{i=1}^n\lambda_i-\nu)^2\\
&\mbox{s.t.}&\lambda_i\geq 0,\ i = 1,...,n.
\end{array}
\end{eqnarray}
The portion of (P) corresponding to (\ref{dualofdual}) becomes
\begin{eqnarray}\label{SDCDWPconstraint}
\begin{array}{ll}
\min\limits_{w}&\frac{1}{2}\left(
\sum\limits_{i=1}^n \frac{1}{2}w_i^2-\nu\right)^2
+\sum\limits_{i=1}^n \{\frac{\alpha_i}{2}w_i^2
-\psi_iw_i\}\\
\mbox{s.t.}&\psi_iw_i\geq0,\ \
i=1,...,n
\end{array}
\end{eqnarray}
under the nonlinear one-to-one map:
\begin{equation}\label{x-y}
w_i=\left\{
\begin{array}{cl}\sqrt{2\lambda_i},&{\rm if}~\psi_i\geq 0,\\
-\sqrt{2\lambda_i},&{\rm if}~\psi_i< 0,
\end{array}\right\}\  i = 1,...,n.
\end{equation}

From Lemma \ref{lem:4}, we know that the portion specified by
(\ref{SDCDWPconstraint}) contains the global minimizer $w^*$.
However, due to the opposite sign behavior on the first
component, Corollary \ref{cor:3} implies that the local non-global
minimizer $\underline{w}$ is not in that portion. The mapping (\ref{x-y}) was used
to reveal the hidden convexity of (QP1QC) in Part I, but the local
non-global minimizer is definitely excluded from the
transformation. The missing of the local non-global minimizer can been
seen clearly in Examples 1 and 2 of Part I.


%
%
%
%
%
%

\section{Characterization of local maximizer}

It is not difficult to see that the global maximum of problem (\ref{d-well:4}) goes to $+\infty$ as
$\Vert w\Vert^2$ grows without a bound. Hence there is no global
maximizer of the problem. In this section, we provide an analytic study of
the local maximizer of the simplified problem
(\ref{d-well:4}). 
\begin{lem} \label{lem:1n}
If $\overline{w}$ is a local maximizer of (\ref{d-well:4}), then
\begin{eqnarray}
&&\nabla g(\overline{w})=(\frac{1}{2}\Vert \overline{w}\Vert^2-\nu)\overline{w}+D\overline{w}-\psi=0,\label{con:1n}\\
&&\nabla^2 g(\overline{w})=
\overline{w}\overline{w}^{T}+(\frac{1}{2}\Vert
\overline{w}\Vert^2-\nu)I +D\preceq 0.\label{con:2n}
\end{eqnarray}
\end{lem}
The proof is easy. Moreover, it follows directly from (\ref{con:2n}) that
\begin{equation}\label{a2}
\frac{1}{2}\Vert \overline{w}\Vert^2-\nu +\alpha_i\leq 0,~i=1,2,...,n.
\end{equation}
In other words, at the local maximizer, the value of the Cauchy-Green strain is
smaller than the negative value of all material constants.

\begin{lem}\label{max:re}
If $\overline{w}$ is a local maximizer of (\ref{d-well:4}), then
\[
\psi_i=0 ~\ if \ and \ only \ if \  \overline{w}_i=0,~i=1,2,\ldots,n.
\]
\end{lem}
\proof
It follows from (\ref{con:1n}) that
\[
(\Vert \overline{w}\Vert^2-2\nu+2\alpha_i)\overline{w}_i=\psi_i.
\]
If $\overline{w}_i=0$, then $\psi_i=0$. On the other hand, if
$\psi_i=0$, it implies from $ (\Vert
\overline{w}\Vert^2-2\nu+2\alpha_i)\overline{w}_i=0$ that either
$\Vert \overline{w}\Vert^2-2\nu+2\alpha_i=0$ or $\overline{w}_i=0$
(or both). Suppose that $\Vert \overline{w}\Vert^2-2\nu+2\alpha_i=0$. It
follows from (\ref{con:2n}) that
$$ \overline{w}_i^{2}+\frac{1}{2}\Vert \overline{w}\Vert^2-\nu
+\alpha_i=\overline{w}_i^{2}\leq 0.$$
Therefore, $\overline{w}_i$ must
be also $0$, and the proof follows.
\endproof

\begin{lem}\label{mc:1}
If $\nu - \alpha_n\le 0$, then the double well potential problem
(\ref{d-well:4}) has no local maximizer.
\end{lem}
\proof If $\nu - \alpha_n<0$, then (\ref{a2}) cannot be true and we have the conclusion.
Now, assume that $\nu - \alpha_n=0$. If (\ref{d-well:4}) has a local
maximizer $\overline{w}$, then it follows from (\ref{a2})
that $\|\overline{w}\|=0$ or, equivalently, $ \overline{w}=0. $
By Lemma \ref{max:re}, we have $\psi=0$. It is routine to verify that
\begin{eqnarray}
&& \frac{\partial g(w)}{\partial w_n}=
(\frac{1}{2}\Vert w\Vert^2-\nu)w_n+\alpha_nw_n-\psi_n, \nonumber\\
&& \frac{\partial^2 g(w)}{\partial^2 w_n}= w_n^2 +
\frac{1}{2}\Vert w\Vert^2-\nu +\alpha_n,\nonumber\\
&& \frac{\partial^3 g(w)}{\partial^3 w_n }= 3w_n.\nonumber
\end{eqnarray}
Since $\overline{w}=0$ and $\psi=0$, we have
\begin{equation}
\frac{\partial g(\overline{w})}{\partial w_n}=\frac{\partial^2
g(\overline{w})}{\partial^2 w_n}=\frac{\partial^3
g(\overline{w})}{\partial^3 w_n}=0.\label{gw_n}
\end{equation}
Notice that
\[
\frac{\partial^4 g(\overline{w})}{\partial^4 w_n }= 3>0.
\]
Consequently, $\overline{w}=0$ is not a local maximizer and we reached a
contradiction. This completes the proof.
\endproof

\begin{lem}\label{mc:2}
If $\nu - \alpha_n> 0$ and $\psi=0$, then the double well potential
problem (\ref{d-well:4}) has a unique local maximizer
$\overline{w}=0$.
\end{lem}
\proof 
Since $\nu - \alpha_n> 0,~\psi=0$ and
$\alpha_1\le\ldots\le\alpha_n$, we have
\begin{eqnarray*}
&&\nabla g(0)=0,\\
&&\nabla^2 g(0)=  -\nu I +D \prec 0.
\end{eqnarray*}
Therefore, $\overline{w}=0$ is a local maximizer of problem
(\ref{d-well:4}). Lemma \ref{max:re} further guarantees that $\overline{w}=0$
is the unique local maximizer.
\endproof

\begin{lem}\label{mc:3}
If $\nu - \alpha_n> 0$ and $\psi\neq 0$, then the double well
potential problem (\ref{d-well:4}) has at most one local maximizer.
\end{lem}
\proof  Suppose $\overline{w}$ is a local maximizer of
(\ref{d-well:4}). Since $\psi\neq 0$, we let $k\in [1:n]$ be the
largest nonzero index in $\{1,...,n\}$ such that
\begin{equation}
\psi_{k}\neq 0;~ \psi_{k+1}=\ldots=\psi_n=0.\label{nonzero}
\end{equation}
In addition, let $I_k$ be the identity matrix of order $k$ and
$D_k={\rm Diag}(\alpha_1,\ldots,\alpha_k).$

From equation (\ref{con:1n}), we have $(\Vert
\overline{w}\Vert^2-2\nu+2\alpha_k)\overline{w}_k=2\psi_k$. Since
$\psi_{k}\neq 0$, by Lemma \ref{max:re}, we know
$\overline{w}_k\neq 0$, which implies that $\Vert
\overline{w}\Vert^2-2\nu +2\alpha_k\not=0$. From inequality (\ref{a2}), we
further know that $\Vert \overline{w}\Vert^2-2\nu +2\alpha_k<0.$ Moreover, we have
\begin{equation}
\Vert \overline{w}\Vert^2-2\nu +2\alpha_i<0,~i=1,...,k.\label{a5}
\end{equation}
Consequently, matrix $(\Vert \overline{w}\Vert^2-2\nu)I_k+2D_k$ is negative
definite and, once $\Vert \overline{w}\Vert$ is computed, $\overline{w}_{i},~i=1,2,\ldots,k,$ can be uniquely
determined by the following system of equations:
\begin{equation} \overline{w}_i=\frac{2\psi_i}{\Vert
\overline{w}\Vert^2-2\nu+2\alpha_i}, ~~i=1,...,k.\label{sol:n}
\end{equation}
 Since
$\psi_{k+1}=\ldots=\psi_n=0$ implies that
$\overline{w}_{k+1}=\ldots=\overline{w}_n=0$, it follows that any
local maximizer $\overline{w}$ must satisfy that
$$
\sum_{i=1}^k\frac{4\psi_i^2}{(\Vert
\overline{w}\Vert^2-2\nu+2\alpha_i)^2} =\Vert \overline{w}\Vert^2.
$$



From (\ref{a3}) and (\ref{a5}), we know  $\Vert \overline{w}\Vert^2$ is a root
of the following convex secular function:
\begin{equation}
 h(t)=\sum_{i=1}^k\frac{4\psi_i^2}{(t-2\nu+2\alpha_i)^2}-t,
 ~~t\in[0,2\nu-2\alpha_k).\label{equ:n}
\end{equation}
Since matrix $(\Vert \overline{w}\Vert^2-2\nu)I_k+2D_k$ is negative
definite, from (\ref{con:2n}), we know that
$$
-2(\Gamma_k \overline{w}^k) (\Gamma_k \overline{w}^k)^T+ I_k
\succeq 0,
$$
where $\overline{w}^k=(\overline{w}_1,\ldots,\overline{w}_k)^T$ and
$$
\Gamma_k={\rm Diag}\left( \frac{1}{\sqrt{-\Vert
\overline{w}\Vert^2+2\nu-2\alpha_1}}, \ldots,\frac{1}{\sqrt{-\Vert
\overline{w}\Vert^2+2\nu-2\alpha_k}} \right).
$$
Then, 
\begin{eqnarray}
0&\leq& {\rm det}\left(-2(\Gamma_k \overline{w}^k) (\Gamma_k\overline{w}^k)^T+ I_k\right)\label{determinant3} \\
&=&-2(\Gamma_k \overline{w}^k)^T(\Gamma_k\overline{w}^k)+1\nonumber \\
&=&\sum_{i=1}^k\frac{8\psi_i^2}{(\Vert
\overline{w}\Vert^2-2\nu+\alpha_i)^3}+1\nonumber\\
&=&-h'(t)\mid_{t=\Vert \overline{w}\Vert^2}.~~ (by\ (\ref{a4}))\label{determinant4}
\end{eqnarray}
Since a strictly convex function can have at most one root with its
first derivative being non-positive, based on (\ref{sol:n}), we can conclude
that there is at most one local maximizer.
\endproof


Combining Lemmas \ref{mc:1}, \ref{mc:2} and \ref{mc:3} together, we
have the next result.
\begin{thm}\label{theorem2}
The double well potential problem (\ref{d-well:4}) has at most one
local maximizer.
\end{thm}

The above result can be further extended to obtain a necessary
and sufficient condition under which a local maximum exists.

\begin{thm}\label{cor:00}
The double well potential problem (\ref{d-well:4}) has a local
maximizer if and only if $\nu-\alpha_n>0$ and there is a
$\overline{t}^*\in[0,2\nu-2\alpha_n)$ such that
 $h(\overline{t}^*)=0$ and $ h'(\overline{t}^*)<0$ for the secular function defined in (\ref{equ:n}).
Moreover, if it exists, the local  maximizer $\overline{w}$ is given by
\begin{equation}
\overline{w}
=\left(\frac{2\psi_1}{\overline{t}^*-2\nu+2\alpha_1},\ldots,\frac{2\psi_n}
{\overline{t}^*-2\nu+2\alpha_n} \right). \label{tw1}
\end{equation}
\end{thm}
\proof
(i) (if part) When $\psi=0$, Lemma \ref{mc:2}
assures that $\overline{w}=0$ is the unique local maximizer of
(\ref{d-well:4}), which can be expressed as (\ref{tw1}).


Now, consider $\psi\neq 0$. Let $k=1,...,n$ be defined as in
(\ref{nonzero}) and $\overline{w}$ as in (\ref{tw1}). Since
$h(\overline{t}^*)=0$ and
$\overline{w}_{k+1}=\ldots=\overline{w}_n=0$, we have
$$\overline{t}^*=\sum_{i=1}^k\frac{4\psi_i^2}{(\overline{t}^*-2\nu+2\alpha_i)^2}
=\Vert\overline{w}\Vert^2\in[0,2\nu-2\alpha_n).$$ Then we see that
$\overline{w}$ satisfies the first order necessary condition
(\ref{con:1n}). Moreover, we have
\begin{equation}
(\frac{1}{2}\Vert \overline{w}\Vert^2-\nu)I+D\prec0.\label{ne}
\end{equation}
Let $\overline{w}^k=(\overline{w}_1,\ldots,\overline{w}_k)^T$. Using
Weyl's inequality (see \cite{Horn}, Theorem 4.3.1), we have
\begin{eqnarray}
\sigma_i\left(\overline{w}^k (\overline{w}^k)^{T}+(\frac{1}{2}\Vert
\overline{w}^k \Vert^2-\nu)I_k+D_k\right)&\le&
\sigma_i\left(\overline{w}^k
(\overline{w}^k)^{T}\right)+\sigma_k\left((\frac{1}{2}\Vert
\overline{w}\Vert^2-\nu)I_k+D_k\right)\nonumber\\
&<&\sigma_i\left(\overline{w}^k (\overline{w}^k)^{T}\right)\nonumber\\
&=&0,~~~~~\hbox{for}~i=1,2,\ldots,k-1.\nonumber
\end{eqnarray}
Therefore, the first $k-1$ eigenvalues of the matrix
$\overline{w}^k(\overline{w}^k)^{T}+(\frac{1}{2}\Vert \overline{w}
\Vert^2-\nu)I_k+D_k$ are negative. It follows from
(\ref{determinant3}), (\ref{determinant4}) and the assumption of
$h'(\overline{t}^*)<0$ that
\begin{equation}
\det\left(-\overline{w}^k(\overline{w}^k)^{T}-(\frac{1}{2}\Vert
\overline{w}\Vert^2-\nu)I_k-D_k\right)=\frac{-
h'(\overline{t}^*)}{2{\det}^2(\Gamma_k)}>0.\label{det:41}
\end{equation}
If $k$ is even, then
$$\det\left(\overline{w}^k(\overline{w}^k)^{T}+(\frac{1}{2}\Vert
\overline{w}\Vert^2-\nu)I_k+D_k\right)=\det\left(-\overline{w}^k
(\overline{w}^k)^{T}-(\frac{1}{2}\Vert \overline{w}
\Vert^2-\nu)I_k-D_k\right)>0,$$ which implies that the $k^{th}$
eigenvalue of matrix $\overline{w}^k(\overline{w}^k)^{T}+(\frac{1}{2}\Vert
\overline{w} \Vert^2-\nu)I_k+D_k$ is negative.

If $k$ is odd, then
$$\det\left(\overline{w}^k (\overline{w}^k)^{T}+(\frac{1}{2}\Vert
\overline{w} \Vert^2-\nu)I_k+D_k\right)=-\det\left(-\overline{w}^k
(\overline{w}^k)^{T}-(\frac{1}{2}\Vert \overline{w}
\Vert^2-\nu)I_k-D_k\right)<0,$$ which says that the $k^{th}$ eigenvalue
of matrix $\overline{w}^k(\overline{w}^k)^{T}+(\frac{1}{2}\Vert
\overline{w} \Vert^2-\nu)I_k+D_k$ is again negative. In other words,
$\overline{w}^k (\overline{w}^k)^{T}+(\frac{1}{2}\Vert \overline{w}
\Vert^2-\nu)I_k+D_k\prec 0$. From (\ref{ne}), we have
\[
\overline{w}
 \overline{w}^{T}+(\frac{1}{2}\Vert \overline{w}
\Vert^2-\nu)I+D=\left[\begin{array}{cc}\overline{w}^k
(\overline{w}^k)^{T}+(\frac{1}{2}\Vert \overline{w}
\Vert^2-\nu)I_k+D_k&0\\
0& (\frac{1}{2}\Vert \overline{w}
\Vert^2-\nu)I_{n-k}+D_{n-k}\end{array}\right]\prec 0,
\]
where $D_{n-k}={\rm Diag}(\alpha_{k+1},\ldots,\alpha_{n})$.
Consequently, $\overline{w}$ satisfies the second order sufficient
condition and becomes a local maximizer.

(ii) (only if part)
Let $\overline{w}$ be the unique local maximizer of
(\ref{d-well:4}). By Lemma \ref{mc:1}, we have $\nu-\alpha_n>0$. If
$\psi=0$, Lemma \ref{mc:2} implies that $\overline{w}=0$ . In
this case, since $h(t)=-t$, there is a unique $\overline{t}^*=0$
such that $h(\overline{t}^*)=0, \overline{t}^*\in[0,2\nu-2\alpha_k)$
and $h'(\overline{t}^*)=-1<0$. The expression (\ref{tw1})
follows immediately.

Assume that $\psi\neq 0$ with $\psi_{k}\neq 0$ and
$\psi_{k+1}=\ldots=\psi_n=0$, and let
$\overline{t}^*=\Vert\overline{w}\Vert^2$. From
(\ref{equ:n}) and (\ref{determinant4}), we know $h(\overline{t}^*)=0,~
\overline{t}^*\in[0,2\nu-2\alpha_k)$ and $h'(\overline{t}^*)\le0$.
Then the expression for $\overline{w}$ in (\ref{tw1}) follows
from (\ref{sol:n}) and $\overline{w}_{k+1}=\ldots=\overline{w}_n=0$.
In the rest of the proof, we shall show a stronger result
of which $\overline{t}^*\in[0,2\nu-2\alpha_n)$ and $
h'(\overline{t}^*)<0.$

First, if $k<n$, then $\psi_n=\overline{w}_n=0$. From
(\ref{a2}), we know $\overline{t}^*\in[0,2\nu-2\alpha_n]$. Suppose that $\overline{t}^*=2\nu-2\alpha_n$, similar to
(\ref{gw_n}), we can verify that
\[
\frac{\partial g(\overline{w})}{\partial w_n}=\frac{\partial^2
g(\overline{w})}{\partial^2 w_n}=\frac{\partial^3
g(\overline{w})}{\partial^3 w_n}=0
\]
and $\frac{\partial^4 g(\overline{w})}{\partial^4 w_n}=3>0.$ This is
a contradiction to the fact that $\overline{w}$ being a local maximizer of probem
(\ref{d-well:4}). Therefore, $\overline{t}^*\in[0,2\nu-2\alpha_n)$.

Next, we show that $ h'(\overline{t}^*)<0.$,
If not so, we consider $ h'(\overline{t}^*)=0.$ By (\ref{det:41}), we have
$$
\det\left(-\overline{w}^k(\overline{w}^k)^{T}-(\frac{1}{2}\Vert
\overline{w}\Vert^2-\nu)I_k-D_k\right)=-
h'(\overline{t}^*)/{\det}^2(\Gamma_k)=0.
$$
Consequently,  matrix $\overline{w}^k(\overline{w}^k)^{T}+(\frac{1}{2}\Vert
\overline{w}\Vert^2-\nu)I_k+D_k $ is singular and there exists a
$u=(u_1,\ldots,u_k)^T\neq 0$ such that
\begin{equation}
\overline{w}^k(\overline{w}^k)^{T}u+\left((\frac{1}{2}\Vert
\overline{w}^k\Vert^2-\nu)I_k+D_k \right)u=0.\label{second-0}
\end{equation}
Equivalently,
\[
u_i=\frac{-\overline{w}_i(u^T\overline{w}^k)}{\frac{1}{2}\Vert
\overline{w}\Vert^2-\nu+\alpha_i},~i=1,2,\ldots,k.
\]
Since $u\neq 0$, we know $u^T\overline{w}^k\neq 0.$ Define
\[
\tilde u=(u_1,\ldots,u_k,0,\ldots,0)^T\in R^n.
\]
Then, we have
\begin{equation}
{\tilde u}^T{\tilde u}\neq 0,~~ {\tilde u}^T\overline{w}=u^T\overline{w}^k\neq
0.\label{nzero:1}
\end{equation}
Similar to the proof of Theorem \ref{cor:0}, we can consider
$q(\beta):=g(\overline{w}+\beta {\tilde u})$. It follows that
\begin{eqnarray*}
q'(\beta)&=& \nabla g(\overline{w}+\beta {\tilde u})){\tilde u}, \\
q''(\beta)&=& {\tilde u}^T \nabla^2g(\overline{w}+\beta {\tilde u}){\tilde u},\\
q'''(\beta)&=&3{\tilde u}^T (\overline{w}+\beta {\tilde u} ){\tilde u}^T{\tilde u}.
\end{eqnarray*}
Since $\overline{w}$ satisfies the first order necessary condition
(\ref{con:1n}), $q'(0)=0$. By (\ref{second-0}), we have $q''(0)=0$. Moreover,
(\ref{nzero:1}) implies that $ q'''(0)=3({\tilde u}^T
\overline{w})({\tilde u}^T{\tilde u})\neq 0. $ Consequently, $0$ is not a
local maximizer of $q(\beta)$ and $\overline{w}$ is not a local
maximizer, which causes a contradiction. Therefore, we know $h'(\overline{t}^*)<0$. This completes the proof.
\endproof

The next result shows that, when it exists, the unique local
maximizer is ``surrounded'' by all local (non-global and global)
minimizers.
\begin{thm}\label{theorem4}
If $\underline{w}$ is a local minimizer and $\overline{w}$ is the
local maximizer of the double well potential problem
(\ref{d-well:4}), then
\begin{equation}
\|\overline{w}\|<\|\underline{w}\|.\label{relat}
\end{equation}
\end{thm}
\proof (i) $n\geq 2$:  If $\alpha_1<\alpha_2$, following Lemma
\ref{lem:2} and (\ref{a2}), we have
\begin{equation*}
\Vert \underline{w}\Vert^2> 2\nu-2\alpha_2\geq \Vert
\overline{w}\Vert^2.
\end{equation*}
Otherwise, $\alpha_1=\alpha_2$ and we assume that
$\|\overline{w}\|=\|\underline{w}\|$. Applying Lemma
\ref{lem:2} and (\ref{a2}) again, we have 
\begin{equation*}
\Vert \underline{w}\Vert^2= 2\nu-2\alpha_1=2\nu-2\alpha_2=
\Vert\overline{w}\Vert^2.
\end{equation*}
Since both $\underline{w}$ and $\overline{w}$ are critical points of
$g(w)$, Lemma \ref{lem:3} implies that both of them are global
minimizers, which is impossible. Therefore, $\|\overline{w}\|<
\|\underline{w}\|$.

(ii) $n=1$: If $\psi_1=0$, then the first order
necessary condition (\ref{con:1n}) implies that either
$\overline{w}=0$ or $\overline{w}^2-2\nu+2\alpha_1=0$. Since
$\overline{w}$ can not be a global minimizer, the latter case is
eliminated and thus $\overline{w}=0$. To prove (\ref{relat}), it is
sufficient to show that $\underline{w}\neq 0$.
Suppose that $\underline{w}=0$, then the second
order necessary condition (\ref{con:2}) implies that
\begin{equation}
0\leq
2\underline{w}^2+\underline{w}^2-2\nu+2\alpha_1=-2\nu+2\alpha_1.
\end{equation}
By Corollary \ref{cor:1} (i), $g(w)$ is convex  and
hence the local maximizer $\overline{w}$ does not exist, which causes a
contradiction to the setting of the theorem.

If $\psi_1\neq 0$, then $w^2-2\nu+2\alpha_1\neq0$
for any local minimizer or maximizer $w$. Therefore,
$t_1=\underline{w}^2$ and $t_2=\overline{w}^2$ are two solutions to the following equation:
\begin{equation*}
 h(t)= \frac{4\psi_1^2}{(t-2\nu+2\alpha_1)^2}-t=0, 
\end{equation*}
From the proofs of Theorem \ref{theorem1} and Theorem
\ref{theorem2}, we have
\begin{eqnarray*}
h'(t)\mid_{t=\underline{w}^2}&\geq& 0,\\
h'(t)\mid_{t=\overline{w}^2}&\leq &0.
\end{eqnarray*}
Since $h(t)$ is strictly convex, it has two distinct solutions
satisfying the above first order conditions only when
$\underline{w}^2> \overline{w}^2.$ This completes the proof.
\endproof

\section{Computational algorithms}

According to Corollaries \ref{cor:0} and \ref{cor:00}, the
local, non-global minimizer and the local maximizer of the simplified
version of (\ref{d-well:4}), if they exist, are closely related to the
convex secular function $h(t)$ over different intervals. The convex
secular function $h(t)$ is a convenient substitute for the first order
necessary condition, while the intervals capturing the root of
$h(t)$ reflect the second order necessary condition. The sign of
the first derivative of $h(t)$ at the root provides necessary and
sufficient conditions for the type of a local extremum, namely, positive
sign for the local, but non-global, minimizer; negative sign for the local
maximizer.

The necessary and sufficient condition for the global minimum $w^*$
in Theorem \ref{thm:suf-nec} can also be expressed in terms of the
secular function $h(t)$. From (\ref{glob:s1}), we have $\|w^*\|^2\in
[2\nu-2\alpha_1,\infty)$. If $\|w^*\|^2>2\nu-2\alpha_1$,
(\ref{glob:s2}) implies that $h(\|w^*\|^2)=0$. Moreover, by
(\ref{a4}),
\[
h'(\|w\|^2) =
-\sum_{i=1}^n\frac{8\psi_i^2}{(\|w\|^2-2\nu+2\alpha_i)^3}-1<0
\hbox{ for }~\|w\|^2>2\nu-2\alpha_1.
\]
It implies that $h(t)$ is monotonically decreasing on
$(2\nu-2\alpha_1,\infty)$ and the unique root $\|w^*\|^2$ must
recover $w^*$.

Otherwise, if $\|w^*\|^2=2\nu-2\alpha_1$, the secular function
$h(t)$ is singular at $\|w^*\|^2$ and there could be multiple global
minimum solutions. In this case, let ${\bar k}$ be the index such
that $\alpha_1=\alpha_2=\ldots=\alpha_{\bar k}<\alpha_{{\bar k}+1}$.
The first order necessary condition $(\frac{1}{2}\Vert
w\Vert^2-\nu+D)w=\psi$ can be solved by letting
$w=(-\alpha_1I+D)^{+} \psi +\sum_{i=1}^{\bar k} \gamma_i e_i$, where
$(\cdot)^+$ denotes the Moore-Penrose generalized inverse;
$\gamma_i$ are free parameters and $e_i$ is the $i$-th column of
$I$. Then, we can establish the following {\it generalized secular
equation}
\begin{equation}\label{secular-eq}
\|w\|^2= \|(- \alpha_1I+D)^{+} \psi +\sum_{i=1}^k
\gamma_i e_i\|^2 =  2\nu - 2\alpha_1
\end{equation}
from which we try to find solution(s)
$\gamma=\gamma^*=(\gamma^*_1,\ldots,\gamma^*_{\bar k})$. Since the
vector $(- \alpha_1I+D)^{+}\psi$ is perpendicular to each vector of
$\gamma^*_i e_i$, we have
\begin{equation}\label{secular-eq:root}
 \sum_{i=1}^k\gamma_i^{*2}={2\nu-2\alpha_1-\| (- \alpha_1I+D)^{+} \psi\|^2}.
\end{equation}
If $\bar k=1$ and ${2\nu-2\alpha_1-\| (- \alpha_1I+D)^{+}
\psi\|^2}>0$, there are exactly two global optimal solutions. If
$\bar k\ge2$ and ${2\nu-2\alpha_1-\| (- \alpha_1I+D)^{+}
\psi\|^2}>0$, there are infinitely many global solutions which form
a $k$-dimensional sphere. The result coincides with Theorem 1 of
Part I. If ${2\nu-2\alpha_1-\| (- \alpha_1I+D)^{+} \psi\|^2}=0$, the
optimal solution set degenerates to a singleton $w^*=(-
\alpha_iI+D)^{+} \psi$.


In summary, we provide three algorithms for finding the global
minimizers, local non-global minimizer and local maximizer,
respectively.

\noindent
\begin{algo} \label{alg1} ~ (finding global minimizers)  \\
\begin{itemize}
\item[Step 1:] Solve the equation of one variable
$$
h(t)=\|[(\frac{1}{2}t-\nu)I+D]^{-1} \psi \|^2-t=0,~ t\in
(2\nu-2\alpha_1,\infty).
$$
If there is a solution $t^{*}$, Stop! The unique global minimizer of
the double well potential problem (\ref{d-well:4}) is
$$
w^{*}=[(\frac{1}{2} t^{*}-\nu)I+D]^{-1}\psi.
$$
Otherwise, go to Step 2.
\item[Step 2]
If $\alpha_1<\alpha_2$ and $\bar k$=1, solve equation (\ref{secular-eq:root}) for
at most two solutions:
$$\gamma_1^{*}=\pm\sqrt{2\nu-2\alpha_1-\| (- \alpha_1I+D)^{+}
\psi\|^2}.$$
If $\gamma_1^{*}\not=0$, the double well potential problem (\ref{d-well:4}) has exactly
two global minimizers of the form
$$
w^{*}=(-\alpha_1 I+D)^{+}\psi+\gamma_1^* e_1.
$$
If $\gamma_1^{*}=0$, $w^{*}=(-\alpha_1 I+D)^{+}\psi$ is the unique global minimizer.
\item[Step 3]
If $\bar k\ge 2$, the double well potential problem
(\ref{d-well:4}) has one or infinitely many global minimizers:
$$
w_i^{*}=(-\alpha_1 I+D)^{+}\psi+\sum_{i=1}^{\bar k} \gamma_i^* e_i,
$$
where $(\gamma^*_1,\ldots,\gamma^*_{\bar k})$ are obtained by solving
(\ref{secular-eq:root}).

If $ \sqrt{2\nu-2\alpha_1-\| (-
\alpha_1I+D)^{+} \psi \|^2}=0,$ $w^{*}=(-\alpha_1 I+D)^{+}\psi$
is the unique optimal solution.

Otherwise, the global optimal
solutions form a sphere centered at $(- \alpha_1I+D)^{+} \psi$
with the radius of $\sqrt{2\nu-2\alpha_1-\| (- \alpha_1I+D)^{+}
\psi \|^2}$.
\end{itemize}
\end{algo}

\noindent
\begin{algo}\label{alg2} ~ (finding local non-global minimizer)  \\
Solve the equation
$$
h(t)=\|[(\frac{1}{2}t-\nu)I+D]^{-1} \psi \|^2-t=0,~ t\in
(\max\{2\nu-2\alpha_2,0\},2\nu-2\alpha_1).
$$
If there is a solution $\underline{t}^*$ such that
$h'(\underline{t}^*)> 0$, the unique  local non-global minimizer of
the double well potential problem (\ref{d-well:4}) is
$$
\underline{w}=[(\frac{1}{2}\underline{t}^*-\nu)I+D]^{-1}\psi.
$$
Otherwise, declare that there is no local non-global minimizer.

\end{algo}

\noindent
\begin{algo}\label{alg3}  ~ (finding the local maximizer) \\
\begin{itemize}
\item[Step 1]
If $ \nu - \alpha_n\le0$, declare that there is no local maximizer.

If $ \nu - \alpha_n>0$ and $\psi=0$, then $0$ is the unique local
maximizer.

Otherwise, go to Step 2.

\item[Step 2] 
Solve the equation
$$
h_+(t)= \|[(\frac{1}{2}t -\nu)I+D]^{+} \psi \|^2-t =0,~t\in[0,
2\nu-2\alpha_n).
$$
If there is a solution $\overline{t}^*$ such that
$h'(\underline{t}^*)<0$, then the unique local maximizer of the double
well potential problem (\ref{d-well:4}) is
\[
\overline{w}^*=[(\frac{1}{2}\overline{t}^*-\nu)I+D]^{+} \psi.
\]

Otherwise, declare that there is no local maximizer.



\end{itemize}
\end{algo}

%

Notice that each of the above three algorithms can be done in a polynomial
time since the main computation involved is to solve the secular equation in
one variable. To illustrate their numerically behavior, we use the
same data set of ($A, B, c, d, f$) of the three examples in
Part I of this paper and apply the space reduction in Section 2 to convert
the testing problems into the format of (\ref{d-well:4}).
%
%
%
%
%

\begin{exam}
(Example 1 of Part I:) Let $n=1$ and $\nu = 14, ~ \alpha_1 =-2,~ \psi_1 =-3,$ the double
well potential problem becomes
$$\min \left\{g(w)=
\frac{1}{2}\left(\frac{1}{2}w^2-14\right)^2-w^2+3w\right\}.$$
The corresponding function $g(w)$ is shown in Figure 2.
\end{exam}
In this example, there are one global minimizer, one local non-global
minimizer and one local maximizer. The secular function
\begin{equation}
h(t)=\frac{36}{(t-32)^2}-t \label{m:1}
\end{equation}
is shown in Figure 2. By finding the root of (\ref{m:1}) in
$(2\nu-2\alpha_1,\infty)=(32,\infty)$, Algorithm \ref{alg1} provides a solution
$t^*=33.0438$ and we find the global minimizer $w^*=-5.7484$
with the value of $-47.1089$. For the local non-global
minimizer, we apply Algorithm \ref{alg2} to find the root of
(\ref{m:1}) in $(\max\{2\nu-2\alpha_2,0\},2\nu-2\alpha_1)=(0,32)$.
Algorithm \ref{alg2} returned $\underline{t}^*=30.9210$ with
$h'(\underline{t}^*)=56.3138>0$, which concluded that the
local non-global minimizer is $\underline{w}^*=5.5607$ with the
value of $-13.1725$. As for the local maximizer, by finding
the root of (\ref{m:1}) in $[0,2\nu-2\alpha_n)=[0,32)$, Algorithm
\ref{alg3} returned $\overline{t}^*=0.0352$ with
$h'(\overline{t}^*)=-0.9978<0$. It led to the local maximizer $\overline{w}^*=0.1877$ with the
value of $98.2814$.

Notice that the signs of the two minimizers, $w^*=-5.7484$ and
$\underline{w}^*=5.5607$, are different, which demonstrates
Corollary \ref{cor:3}. The numerical results also showed that
the local maximizer $\overline{w}^*=0.1877$ locates between the two
minimizers, which is claimed by Theorem \ref{theorem4}.

\begin{figure}[!htbp]\label{fig:1}
\begin{center}
\includegraphics[width=0.5\textwidth]{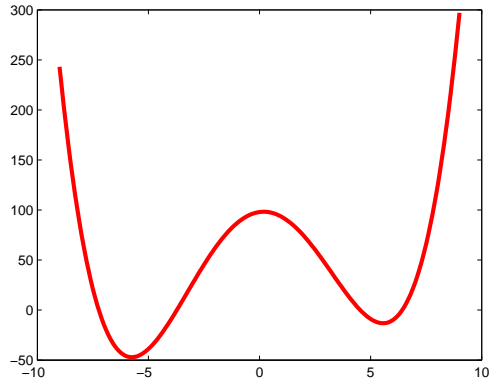}
\end{center}\caption{The graph of $g(w)$ in Example 1 ($n=1$).}
\end{figure}

\begin{figure}[!htbp]\label{fig:1:2}
\begin{center}
\includegraphics[width=0.5\textwidth]{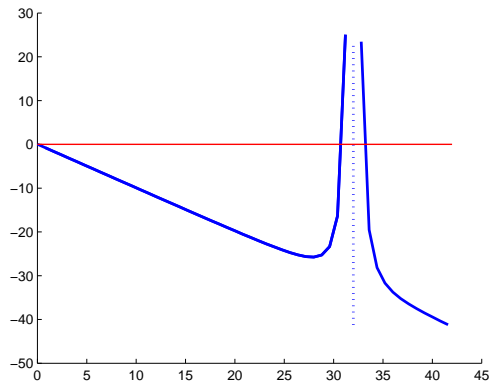}
\end{center}\caption{The secular function (\ref{m:1}).}
\end{figure}

\begin{exam}
(Example 2 of Part I:) Applying the space reduction technique, we
obtain the double well potential problem in the format of
(\ref{d-well:4}) with the data  $n=2$ and
$$
\nu = 27.9994, ~ D
=\left[\begin{array}{cc}-1.9960&0\\0&202.0700\end{array}\right],~
\psi =\left[\begin{array}{c} -22.0487\\-502.0209\end{array}\right].
$$
The corresponding function $g(w)$ and its contour are
shown in Figure 4.
\end{exam}
Its secular function becomes
\begin{equation}
h(t)=\frac{1944.5808}{(t-59.9908)^2}+\frac{1008099.9361}{(t+350.14)^2}-t
\label{m:2}
\end{equation}
(shown in Figure 5). Finding the root of (\ref{m:2}) on
$(2\nu-2\alpha_1,\infty)=(59.9908,\infty)$ results in $t^*=65.6930$,
Algorithm \ref{alg1} gives the global minimizer
   $w^*=\left[\begin{array}{c} -7.7335\\
   -2.4262\end{array}\right]$ with the value of $-841.7182$.
Similarly, finding the root of (\ref{m:2}) in
$(\max\{2\nu-2\alpha_2,0\},2\nu-2\alpha_1)=(0,59.9908)$ results in
$\underline{t}^*= 53.5813$. Since $h'(\underline{t}^*)=13.7390
>0$, Algorithm \ref{alg2} provides the local non-global
minimizer  $\underline{w}^*=\left[\begin{array}{c}  6.8800\\
   -2.4993\end{array}\right]$ with the value of $-518.3996$.
Notice that the signs of the first component of the two minimizers
are different, which demonstrates Corollary \ref{cor:3}. Finally, since
$2\nu-2\alpha_n=-348.1412<0$, 
Algorithm \ref{alg3} says that there is no
local maximizer for this example.

\begin{figure}[!htbp]\label{fig:2}
\begin{center}
\includegraphics[width=0.5\textwidth]{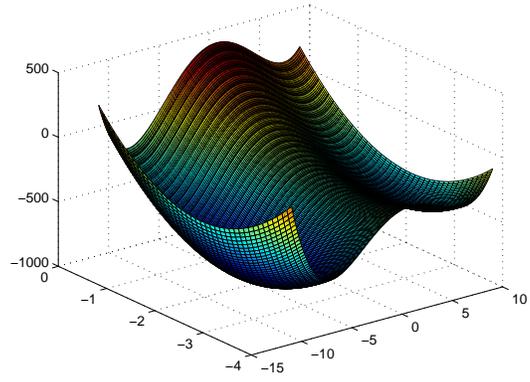}
\includegraphics[width=0.5\textwidth]{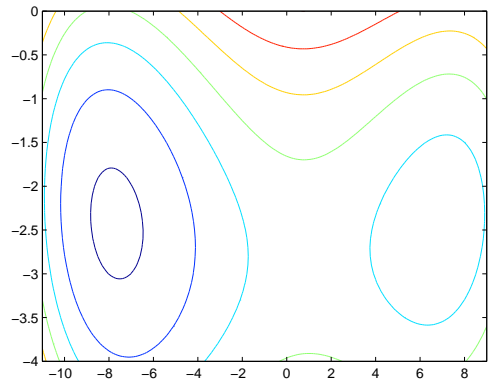}
\end{center}\caption{ The function $g(w)$ in Example 2 and its contour ($n=2$).}
\end{figure}

\begin{figure}[!htbp]\label{fig:2:2}
\begin{center}
\includegraphics[width=0.5\textwidth]{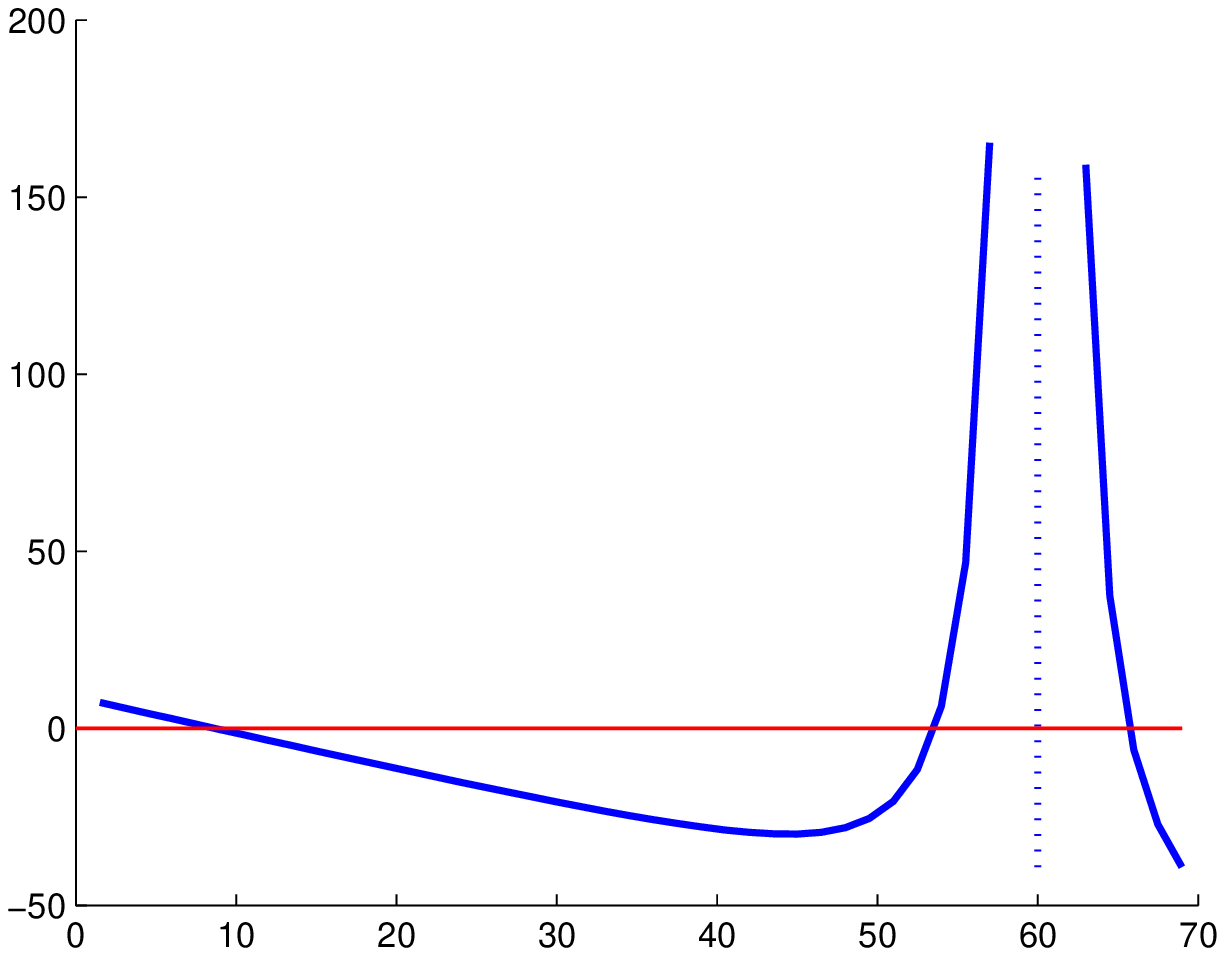}
\end{center}\caption{ The secular function (\ref{m:2}).}
\end{figure}

\begin{exam} (Maxican Hat Example)  In this example,
$$g(w)=\frac{1}{2}(\frac{1}{2}w_1^2+\frac{1}{2}w_2^2-38)^2,$$
which is already in the format of (\ref{d-well:4}) with $n=2$,
$$
\nu = 38, ~ D =\left[\begin{array}{cc}0&0\\0&0\end{array}\right],~
\psi =\left[\begin{array}{c} 0\\0\end{array}\right].
$$
The graph of the Maxican hat function $g(w)$ and its contour are
shown in Figure 6.
\end{exam}

Since $\alpha_1=\alpha_2=0$ and $2\nu-2\alpha_1=76$, the
secular function
\begin{equation}
h(t)=-t,~~t\not=76\label{m:3}
\end{equation}
has a unique solution $0$ and it becomes singular at $t=76$.
Algorithm \ref{alg1} stopped at Step 3 and claimed that
$$
w^{*}=\left\{(\gamma_1^*, \gamma_2^*)\mid ~
(\gamma_1^*)^2+(\gamma_2^*)^2=76\right\}.
$$
is the set of global optimal solutions with the optimal value of $0$.

Since $(\max\{2\nu-2\alpha_2,0\},2\nu-2\alpha_1)=(76,76)=\emptyset$,
Algorithm \ref{alg2} returned an answer that there is no local non-global
minimizer. It is clear that (\ref{m:3}) has a unique
root $\overline{t}^*=0$ on $[0, 2\nu-2\alpha_n)=[0,76)$ and
$h'(\overline{t}^*)= -1<0$. Since $\psi=0$, Algorithm \ref{alg3}
returned
the unique local maximizer $\overline{w}^*=\left[\begin{array}{c} 0\\
0\end{array}\right]$.

\begin{figure}[!htbp]\label{fig:3}
\begin{center}
\includegraphics[width=0.5\textwidth]{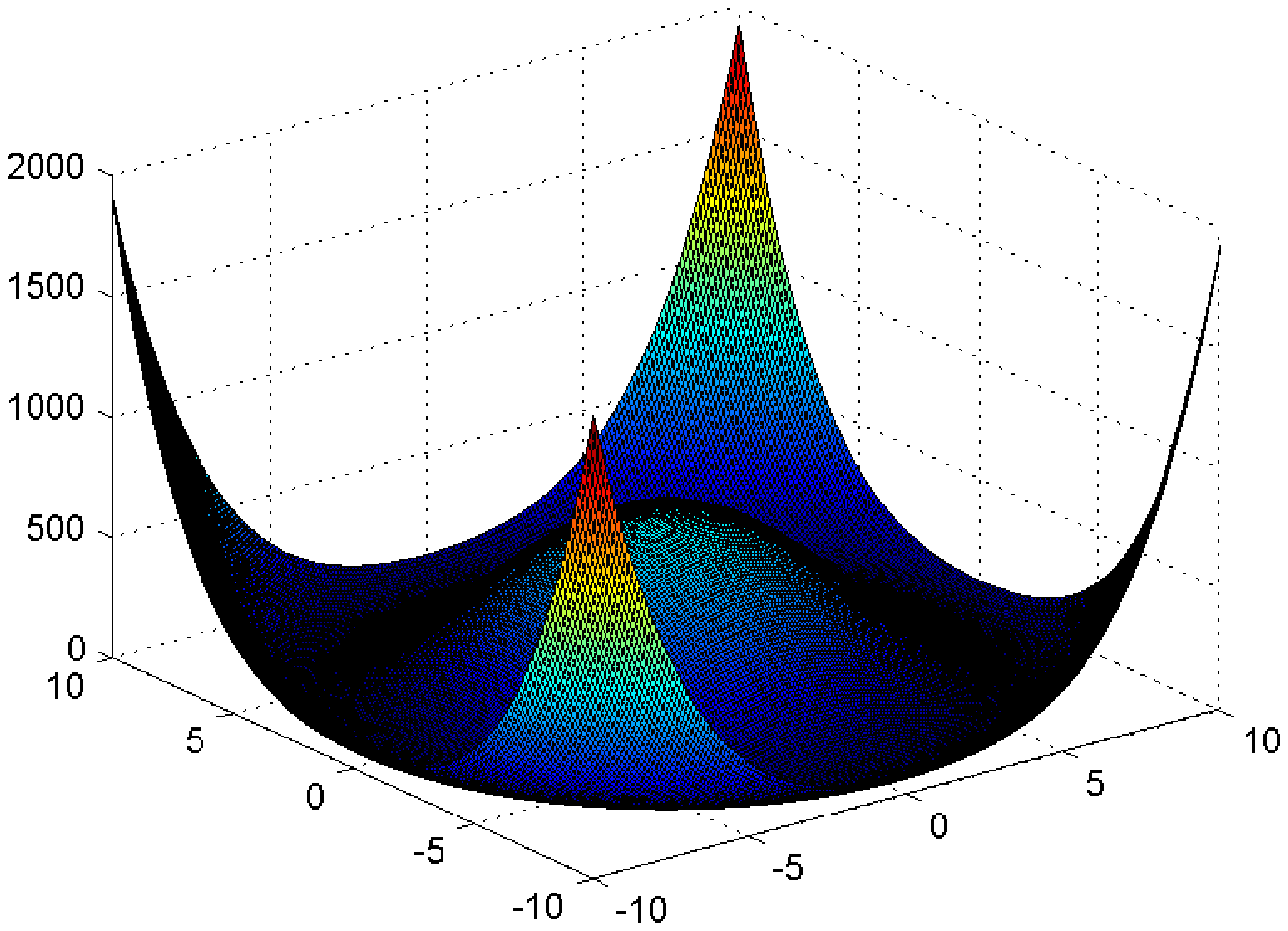}
\includegraphics[width=0.5\textwidth]{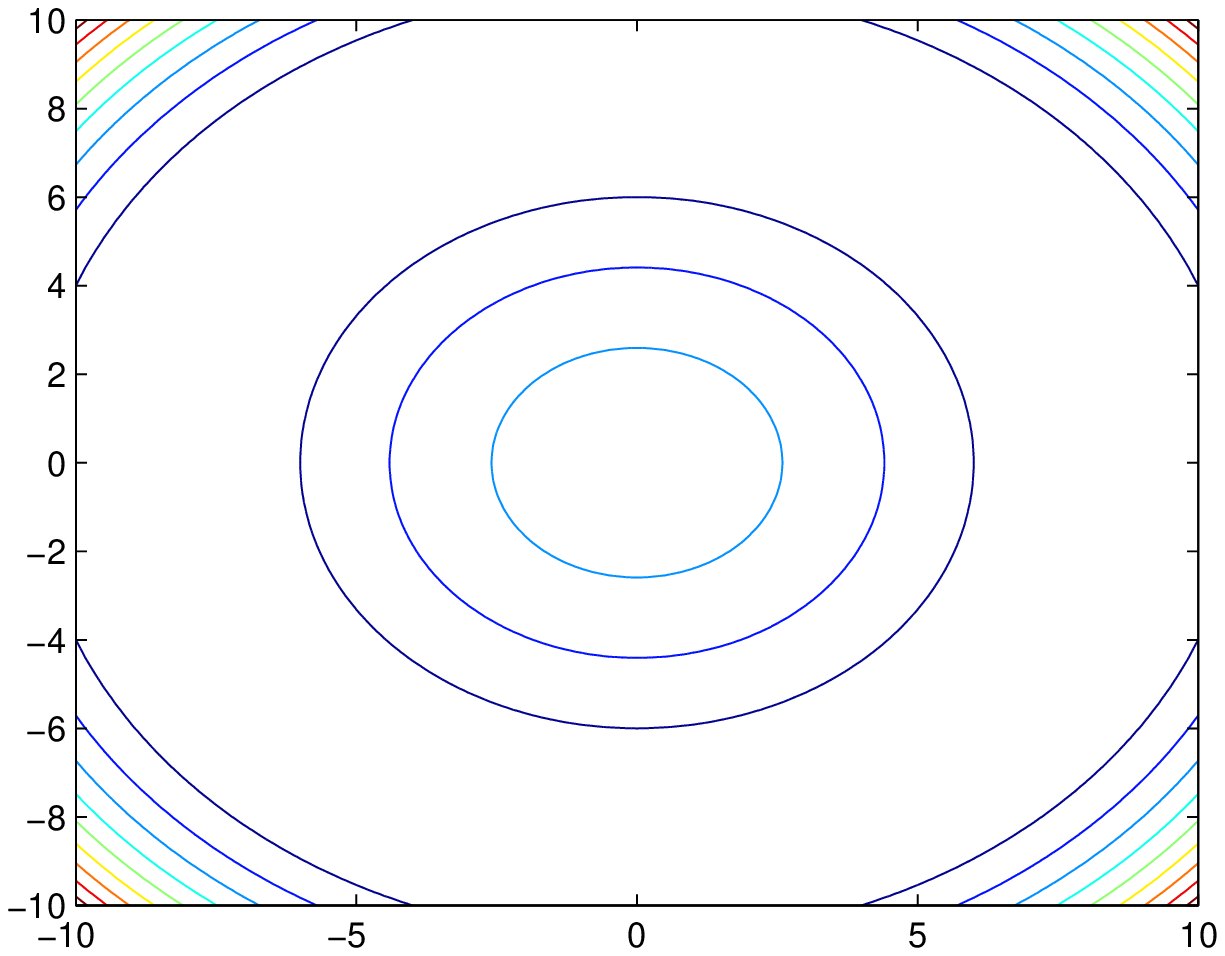}
\end{center}\caption{ The function $g(w)$ in Example 3 and its contour ($n=2$).}
\end{figure}

\section{Conclusions}

In this paper we have characterized the local minimizers and maximizers of
the double well potential problem. By analyzing the first and the
second order necessary conditions and through the study of the corresponding
secular functions, we are able to estimate the number of local
optimizers and locate each of them. Moreover, the convex secular
functions (equations) are used to characterize sufficient and
necessary conditions for all types of optimizers with explicit computational
 algorithms developed for finding them. The (DWP) problem is a
special case of the more general quadratic programming problem with one
quadratic constraint (QP1QC). We expect that the analytical
techniques developed in this paper can be extended to study (QP1QC) and
 other quadratic
programming problems.

\section*{Acknowledgments}
This research was undertaken while Y. Xia visited the National Cheng
Kung University, Tainan, Taiwan. Sheu's research work was sponsored
partially by Taiwan NSC 98-2115-M-006 -010 -MY2 and by National
Center for Theoretic Sciences (The southern branch). Xia's research
work was supported by National Natural Science Foundation of China
under grants 11001006 and 91130019/A011702, and by the fund of State
Key Laboratory of Software Development Environment under grant
SKLSDE-2013ZX-13. Xing's research work was supported by NSFC No.
11171177.

\end{document}